\documentclass[a4paper, 11pt, twoside]{article}
\usepackage{amsmath}
\usepackage{amscd}
\usepackage{amssymb}
\usepackage{amsthm}
\usepackage{graphics}
\usepackage{indentfirst}
\usepackage{latexsym}
\usepackage{amsfonts}
\usepackage{multicol}
\usepackage{stmaryrd}
\usepackage{array}
\usepackage{pxfonts}
\usepackage{youngtab}

\newtheoremstyle{mattthm}{}{}{\itshape}{}{\bfseries}{.}{ }{}

\theoremstyle{mattthm}
\newtheorem{lemma}{Lemma}[section]
\newtheorem{propn}[lemma]{Proposition}
\newtheorem{thm}[lemma]{Theorem}
\newtheorem{cory}[lemma]{Corollary}
\newtheorem{conj}[lemma]{Conjecture}

\newtheoremstyle{mattdef}{}{}{}{}{\bfseries}{.}{ }{}

\theoremstyle{mattdef}
\newtheorem*{rmk}{Remark}

\newtheorem*{eg}{Example}

\newtheorem*{acks}{Acknowledgements}

\begin{document}

\newenvironment{pf}{\noindent\textbf{Proof.}}{\hfill \qedsymbol\newline}
\newenvironment{pfof}[1]{\vspace{\topsep}\noindent\textbf{Proof of {#1}.}}{\hfill \qedsymbol\newline}
\newenvironment{pfenum}{\noindent\textbf{Proof.}\indent\begin{enumerate}\vspace{-\topsep}}{\end{enumerate}\vspace{-\topsep}\hfill \qedsymbol\newline}
\newenvironment{pfnb}{\noindent\textbf{Proof.}}{\newline}

\newlength\raiser
\newcommand\rh[1]{\raisebox{-3pt}{#1}}
\newcommand\mo{-\negthinspace1}
\newcommand\hdom{\linethickness{0.5mm}\put(0,0){\line(1,0){24}}\put(0,12){\line(1,0){24}}\put(0,0){\line(0,1){12}}\put(24,0){\line(0,1){12}}}
\newcommand\vdom{\linethickness{0.5mm}\put(0,0){\line(0,1){24}}\put(12,0){\line(0,1){24}}\put(0,0){\line(1,0){12}}\put(0,24){\line(1,0){12}}}
\newcommand\lra\longrightarrow
\newcommand\rez{\operatorname{res}}
\newcommand\bsm{\begin{smallmatrix}}
\newcommand\esm{\end{smallmatrix}}
\newcommand{\rt}[1]{\rotatebox{90}{$#1$}}
\newcommand\la\lambda
\newcommand{\ol}{\overline}
\newcommand{\ul}{\underline}
\newcommand\reg{^{\operatorname{reg}}}
\newcommand{\lan}{\langle}
\newcommand{\ran}{\rangle}
\newcommand\fkh{\mathfrak{h}}
\newcommand\fka{\mathfrak{a}}
\newcommand\partn{\mathcal{P}}
\newcommand{\py}[3]{\,_{#1}{#2}_{#3}}
\newcommand{\pyy}[5]{\,_{#1}{#2}_{#3}{#4}_{#5}}
\newcommand{\thmlc}[3]{\textup{\textbf{(\!\! #1 \cite[#3]{#2})}}}
\newcommand{\sss}{\mathfrak{S}_}
\newcommand{\dom}{\trianglerighteqslant}
\newcommand{\doms}{\vartriangleright}
\newcommand{\ndom}{\ntrianglerighteqslant}
\newcommand{\ndoms}{\not\vartriangleright}
\newcommand{\domby}{\trianglelefteqslant}
\newcommand{\domsby}{\vartriangleleft}
\newcommand{\ndomby}{\ntrianglelefteqslant}
\newcommand{\ndomsby}{\not\vartriangleleft}
\newcommand{\subs}[1]{\subsection{#1}}
\newcommand{\nin}{\notin}
\newcommand{\nchar}{\operatorname{char}}
\newcommand{\thmcite}[2]{\textup{\textbf{\cite[#2]{#1}}}\ }
\newcommand\zez{\mathbb{Z}/e\mathbb{Z}}
\newcommand\zepz{\mathbb{Z}/(e+1)\mathbb{Z}}
\newcommand{\bbf}{\mathbb{F}}
\newcommand{\bbc}{\mathbb{C}}
\newcommand{\bbn}{\mathbb{N}}
\newcommand{\bbq}{\mathbb{Q}}
\newcommand{\bbz}{\mathbb{Z}}
\newcommand\zo{\bbn_0}
\newcommand{\gs}{\geqslant}
\newcommand{\ls}{\leqslant}
\renewcommand\leq\ls
\renewcommand\geq\gs
\newcommand\dw{^\triangle}
\newcommand\wod{^\triangledown}
\newcommand{\hhh}{\mathcal{H}_}
\newcommand{\hh}{\mathcal{H}}
\newcommand\hsn{\hhh{\bbf,\mo}(\sss n)}
\newcommand\hfq{\hhh{\bbf,q}(\sss n)}
\newcommand{\sect}[1]{\section{#1}}
\newcommand{\ff}{\mathfrak{f}}
\newcommand{\fff}{\mathfrak{F}}
\newcommand\cf{\mathcal{F}}
\newcommand\fkn{\mathfrak{n}}
\newcommand\sx{x}
\newcommand\bra[1]{|#1\ran}
\newcommand\arb[1]{\widehat{\bra{#1}}}
\newcommand\foc[1]{\mathcal{F}_{#1}}
\newcommand{\clam}{\begin{description}\item[\hspace{\leftmargin}Claim.]}
\newcommand{\prof}{\item[\hspace{\leftmargin}Proof.]}
\newcommand{\malc}{\end{description}}
\newcommand\ppmod[1]{\ (\operatorname{mod}\ #1)}
\newcommand\wed\wedge
\newcommand\wede\barwedge
\newcommand\uu[1]{\,\begin{array}{|@{\,}c@{\,}|}\hline #1\\\hline\end{array}\,}
\newcommand{\ux}[1]{\operatorname{ht}_{#1}}
\newcommand\erim{\operatorname{rim}}
\newcommand\mire{\operatorname{rim}'}
\newcommand\mmod{\ \operatorname{Mod}}
\newcommand\call{\mathcal{L}}
\newcommand\calt{\mathcal{T}}
\newcommand\calu{\mathcal{U}}
\newcommand\calv{\mathcal{V}}
\newcommand\calf{\mathcal{F}}
\newcommand\cgs\succcurlyeq
\newcommand\cls\preccurlyeq
\newcommand\inc{\operatorname{asc}}
\newcommand\qbinom[2]{\left[\begin{smallmatrix}#1\\#2\end{smallmatrix}\right]}
\newcommand\dqbinom[2]{\left[\begin{array}c#1\\#2\end{array}\right]}
\newcommand\spe[1]{S^{#1}_{\bbf,\mo}}
\newcommand\smp[1]{D^{#1}_{\bbf,\mo}}
\newcommand\per[1]{M^{#1}_{\bbf,\mo}}
\newcommand\speq[1]{S^{#1}_{\bbf,q}}
\newcommand\smpq[1]{D^{#1}_{\bbf,q}}
\newcommand\perq[1]{M^{#1}_{\bbf,q}}
\newcommand\bx{\put(0,0){\line(0,1){12}}\put(0,12){\line(1,0){12}}\put(12,0){\line(0,1){12}}\put(0,0){\line(1,0){12}}}

\Yboxdim{12pt}

\title{Some reducible Specht modules for\\ Iwahori--Hecke algebras of type $A$ with $q=-\negthinspace1$}
\author{Matthew Fayers\\
\normalsize Queen Mary, University of London, Mile End Road, London E1 4NS, U.K.\\
\\
\large Sin\'{e}ad Lyle\\
\normalsize School of Mathematics, University of East Anglia, Norwich NR4 7TJ, U.K.}\date{}
\maketitle
\begin{center}
\end{center}
\markboth{Matthew Fayers and Sin\' ead Lyle}{Some reducible Specht modules for Iwahori--Hecke algebras of type $A$ with $q=\mo$}
\pagestyle{myheadings}

\begin{abstract}
The reducibility of the Specht modules for the Iwahori--Hecke algebras in type $A$ is still open in the case where the defining parameter $q$ equals $\mo$.  We prove the reducibility of a large class of Specht modules for these algebras.  
\end{abstract}

\sect{Introduction}
Let $n$ be a non-negative integer, $\bbf$ a field, and $q$ an element of $\bbf$.  The \emph{Iwahori--Hecke algebra} $\hh=\hfq$ is a finite-dimensional $\bbf$-algebra which arises in various mathematical contexts.  Its representation theory bears a close resemblance to the representation theory of the group algebra $\bbf\sss n$ (which arises in the special case $q=1$).  A particularly important class of modules for $\hh$ is the class of \emph{Specht modules}; these arise as cell modules for a certain cellular basis of $\hh$, and in cases where $\hh$ is semi-simple the Specht modules are irreducible and afford all the irreducible representations of $\hh$.

In the non-semi-simple case (where $q$ is a root of unity in $\bbf$), it is still interesting to know which Specht modules are irreducible; for the case where $q=1$ and $\bbf$ has prime characteristic $p$, this amounts to asking which ordinary irreducible representations of the symmetric group remain irreducible in characteristic $p$.  The classification of irreducible Specht modules has been studied by several authors, and is almost complete.  This paper is a contribution towards completing the remaining open case, namely the case where $q=\mo\in\bbf$ and the characteristic of $\bbf$ is not $2$.  Our main result is Theorem \ref{main}, where we prove the reducibility of a large class of Specht modules.  
We hope to be able to extend our results in future.

We now give an indication of the layout of this paper.  In Section \ref{mainthmsec}, we give some very basic definitions and state our main result; we also present a conjectured classification of the irreducible Specht modules in the case where $\bbf$ has infinite characteristic.  In Section \ref{backsec}, we recall the additional definitions and background theory that we shall need.  In Section \ref{homsec}, we state some fundamental results on homomorphisms between various modules for $\hh$; we use these to prove further results which aid us in proving reducibility of Specht modules.  In Section \ref{focksec}, we recall the Fock space representation of the quantum group $\calu_v(\widehat{\mathfrak{sl}}_2)$ and its applications to representation theory of Iwahori--Hecke algebras, and we use these results to show reducibility of certain Specht modules.  Finally in Section \ref{pfsec}, we combine our results to complete the proof of the main theorem.

\begin{acks}
The first author was financially supported by a Research Fellowship from the Royal Commission for the Exhibition of 1851; he is very grateful to the Commission.  This work was carried out while the first author was a visiting Postdoctoral Fellow at the Massachusetts Institute of Technology; he is grateful to Richard Stanley for the invitation, and to M.I.T. for its hospitality.

Conjecture \ref{mainconj} in this paper was made while the first author was working with Andrew Mathas at the University of Sydney in 2004; the first author is very grateful for the invitation.

Some of this work was carried out at MSRI Berkeley in March 2008, during the concurrent programmes `Combinatorial representation theory' and `Representation theory of finite groups and related topics'.  Both authors acknowledge generous financial support from MSRI, and wish to thank the organisers of these excellent programmes.
\end{acks}

\sect{The main result and a conjecture}\label{mainthmsec}

In this section, we give our main theorem, and also present a conjectured classification of irreducible Specht modules in the case where $\nchar(\bbf)=\infty$.  First we review the background required to enable us to state our results.

\subsection{Iwahori--Hecke algebras, partitions and Specht modules}

Throughout this paper $\bbf$ denotes a fixed field and $q$ a non-zero element of $\bbf$.  We define $e$ to be the multiplicative order of $q$ in $\bbf$ if $q\neq1$, or the characteristic of $\bbf$ if $q=1$; we adopt the convention that a field whose prime subfield is $\bbq$ has infinite characteristic.  In this paper we shall be primarily concerned with the case where $e=2$ (that is, $q=\mo\in\bbf$), though we shall state results for arbitrary values of $q$ as long as it is convenient.  More general results than those quoted can easily be found elsewhere, especially in the book by Mathas \cite{mathbook}, which is our main reference.  Note, however, that we do not always follow Mathas's conventions; in particular, we use the Specht modules defined by Dipper and James \cite{dija} rather than those of Mathas.

Given any integer $n\gs1$, the \emph{Iwahori--Hecke algebra of the symmetric group $\sss n$} is defined to be the unital associative $\bbf$-algebra $\hfq$ with generators $T_1,\dots,T_{n-1}$ and relations
\begin{alignat*}2
(T_i-q)(T_i+1) &= 0 & &(1\ls i\ls n-1)\\
T_iT_{i+1}T_i &= T_{i+1}T_iT_{i+1} &\qquad&(1\ls i\ls n-2)\\
T_iT_j &= T_jT_i &&(1\ls i<j-1\ls n-2).
\end{alignat*}

The combinatorics describing the representation theory of $\hfq$ is based on compositions and partitions.  A \emph{composition} of $n$ is defined to be a sequence $\la=(\la_1,\la_2,\dots)$ of non-negative integers such that the sum $|\la|=\la_1+\la_2+\dots$ equals $n$; if in addition we have $\la_1\gs\la_2\gs\dots$, we say that $\la$ is a \emph{partition} of $n$.  When writing compositions and partitions, we often omit trailing zeroes and group together equal non-zero parts, and we write $\varnothing$ for the unique partition of $0$.  If $\la$ is a partition, we write $\la'$ for the conjugate partition to $\la$; this is the partition in which
\[\la'_i = \left|\left\{j\in\bbn\ \left|\ \la_j\gs i\right\}\right.\right|.\]
We say that a partition $\la$ is \emph{$e$-regular} if there does not exist $i\gs1$ such that $\la_i=\la_{i+e-1}>0$, and we say that $\la$ is \emph{$e$-restricted} if there is no $i$ with $\la_i-\la_{i+1}\gs e$.

With a partition $[\la]$ is associated its \emph{Young diagram}, which is the set
\[[\la]=\left\{(i,j)\in\bbn^2\ \left|\ j\ls \la_i\right.\right\}.\]
We refer to elements of $\bbn^2$ as \emph{nodes}, and to elements of $[\la]$ as \emph{nodes of $\la$}.

Now we can describe some modules for $\hfq$.  For any composition of $n$, one defines a module $M^\la_{\bbf,q}$ known as the \emph{permutation module}.  If $\la$ is a partition, then $M^\la_{\bbf,q}$ has a distinguished submodule $S^\la_{\bbf,q}$ called the \emph{Specht module}, which is the main object of study in this paper.  We retain the subscript $\bbf,q$ in our notation, to enable us to make statements about modules without reference to the underlying Iwahori--Hecke algebra; for example, when we say that $S^\la_{\bbf,q}$ is irreducible, we mean that it is irreducible as an $\hfq$-module, where $n=|\la|$.

\subsection{The main result}\label{irrsumm}

The purpose of this paper is to consider the question of which Specht modules for $\hfq$ are irreducible.  As with many statements about the representation theory of $\hfq$, the classification of irreducible Specht modules for $\hfq$ involves the parameter $e$ defined above, which is often called the `quantum characteristic'.  If $e=\infty$, then all the Specht modules are irreducible, and they afford all irreducible representations of $\hfq$ as $\la$ varies over the set of partitions of $n$.  So we assume from now on that $e$ is finite.

In the case where $e>2$, the classification of the irreducible Specht modules for $\hfq$ is complete; this result was proved by the authors and others over the course of several papers \cite{jm-3,slred,jm2,mfreduc,mfirred,jlm,slcp}.  To describe the classification, we begin with the case of Specht modules labelled by $e$-regular partitions; the irreducibility of these Specht modules is determined by a theorem known as the \emph{Carter Criterion} \cite[Corollary 5.43]{mathbook}, of which a special case appears in Theorem \ref{Carter} below.  Applying a theorem concerning conjugate partitions (Corollary \ref{conjirr} below), one obtains a corresponding result for Specht modules labelled by $e$-restricted partitions.  The general case is a natural combination of the $e$-regular case and the $e$-restricted case; a partition labelling an irreducible Specht module (called a \emph{JM-partition}) consists of an $e$-regular partition and an $e$-restricted partition, each labelling an irreducible Specht module, joined together in a simple way.

When $e=2$, the Carter Criterion is still valid, so the classification of irreducible Specht modules labelled by $2$-regular or $2$-restricted partitions is known.  However, for the case of partitions which are neither $2$-regular nor $2$-restricted, things are different; one cannot just take the definition of JM-partitions and set $e=2$.  In fact, the partitions of $n$ which are neither $2$-regular nor $2$-restricted and label irreducible Specht modules for $\hsn$ seem to take a very different form from JM-partitions.  Roughly speaking, JM-partitions tend to be `thin', by which we mean that a JM-partition has very few diagonal nodes $(i,i)$ relative to its size; conversely, the partitions labelling irreducible Specht modules when $e=2$ tend to be closer to rectangular partitions (indeed, the rectangular partitions $(a^b)$ all label irreducible Specht modules when $e=2$ and $\nchar(\bbf)=\infty$).  Our main theorem illustrates this difference, since it implies in particular that when $e=2$ a Specht module labelled by a (neither $2$-regular nor $2$-restricted) JM-partition is reducible.

To give our main result, we need to define \emph{ladders}: for $l\gs1$, the $l$th ladder in $\bbn^2$ is defined to be the set
\[\call_l = \left\{\left.(i,j)\in\bbn^2\ \right|\ i+j=l+1\right\}.\]
Given a partition $\la$, we define the $l$th ladder of $\la$ to be the intersection $\call_l(\la)=\call_l\cap[\la]$. We say that $\call_l(\la)$ is \emph{disconnected} if the nodes in $\call_l(\la)$ do not form a consecutive subset of $\call_l$; that is, there exist $1\ls a<b<c\ls l$ such that $(a,l+1-a)$ and $(c,l+1-c)$ lie in $[\la]$, but $(b,l+1-b)$ does not.  Now we can state our main theorem.

\vspace{\topsep}

\noindent\hspace{-3pt}\fbox{\parbox{469pt}{\vspace{-\topsep}

\begin{thm}\label{main}
Suppose $\bbf$ is any field, and $\la$ is a partition.  If there is some $l$ such that the $l$th ladder of $\la$ is disconnected, then the Specht module $\spe\la $ is reducible.
\end{thm}
\vspace{-\topsep}}}
\vspace{\topsep}

The reader may prefer a statement of Theorem \ref{main} that does not involve ladders: it is a simple exercise to show that a partition $\la$ has a disconnected ladder if and only if there exist $1\ls a<b$ such that $\la_a-\la_{a+1}\gs2$ and $\la_b=\la_{b+1}>0$.

Theorem \ref{main} will be proved in the subsequent sections.  For the rest of this section, we consider how to extend it to give a complete classification of irreducible Specht modules.

\subsection{A conjecture in infinite characteristic}

In the case where $\nchar(\bbf)=\infty$, the decomposition numbers for $\hfq$ may be computed using the LLT algorithm \cite{llt}; so for any $\la$, there is a finite algorithm to determine whether $\speq\la$ is reducible.  For the case $e=2$, Andrew Mathas and the first author have carried out these computations for partitions of size at most $45$, and on the basis of this have made a conjecture.

In order to state this conjecture, we need to introduce some more terminology concerning Young diagrams.  If $\la$ is a partition, then we say that a node $(i,j)$ of $\la$ is \emph{removable} if $[\la]\setminus\{(i,j)\}$ is the Young diagram of a partition (i.e.\ if $j=\la_i>\la_{i+1}$), while a node $(i,j)$ not in $[\la]$ is an \emph{addable node of $\la$} if $[\la]\cup\{(i,j)\}$ is the Young diagram of a partition.  For any node $(i,j)$ in $\bbn^2$, we define its \emph{residue} to be the residue modulo $2$ of the integer $j-i$.

Now suppose $\la$ is neither $2$-regular nor $2$-restricted.  Let $a$ be maximal such that $\la_a-\la_{a+1}\gs2$, let $b$ be maximal such that $\la_b=\la_{b+1}>0$, and let $c$ be maximal such that $\la_{a+c}>0$.  We say that $\la$ is an \emph{FM-partition} if the following conditions hold:
\begin{itemize}
\item
$\la_i-\la_{i+1}\ls1$ for all $i\neq a$;
\item
$\la_b\gs a-1\gs b$;
\item
$\la_1>\dots>\la_c$;
\item
if $c=0$, then all the addable nodes of $\la$, except possibly those in the the first row and first column, have the same residue;
\item
if $c>0$, then all addable nodes of $\la$ have the same residue.
\end{itemize}

Now we can give our conjecture; note that the case where a partition is $2$-regular or $2$-restricted is covered by the discussion in \S\ref{irrsumm}, so we can restrict attention to partitions which are neither $2$-regular nor $2$-restricted.

\begin{conj}\label{mainconj}
Suppose $\nchar(\bbf)=\infty$ and that $\la$ is neither $2$-regular nor $2$-restricted.  Then the Specht module $\spe\la$ is irreducible if and only if either $\la$ or $\la'$ is an FM-partition.
\end{conj}

This still leaves open the case where $\la$ has prime characteristic.  Thanks to the theory of decomposition maps (see Theorem \ref{adjmat} below), we know that the set of partitions labelling irreducible Specht modules for $q=\mo$ in characteristic $p$ is a subset of the set of partitions labelling irreducible Specht modules in infinite characteristic.  Experimental evidence suggests that it is a rather small subset; in fact, it seems likely that for any prime $p$ there are only finitely many partitions which are neither $2$-regular nor $2$-restricted and label irreducible Specht modules.  This statement has been proved in the case $p=2$ by James and Mathas \cite{jm2}; here the only such partition is $(2^2)$.  We hope to be able to make a more precise statement in the future.

\sect{Useful background results}\label{backsec}

In this section we summarise some simple background results which we shall need in order to prove our Theorem \ref{main}.

\subsection{Irreducible modules for $\hfq$ and the dominance order}

In order to examine the reducibility of Specht modules, it will be helpful to understand the classification of irreducible $\hfq$-modules.  Let $e$ be as defined in \S\ref{irrsumm}, and suppose that $\la$ is a partition of $n$.  If $\la$ is $e$-regular, then the Specht module $S^\la_{\bbf,q}$ has an irreducible cosocle which is labelled $D^\la_{\bbf,q}$; the modules $D^\la_{\bbf,q}$ give all the irreducible $\hfq$ modules, as $\la$ ranges over the set of $e$-regular partitions of $n$.

The classification of irreducible Specht modules is a special case of the \emph{decomposition number problem}, which asks for the composition multiplicities $[\speq\la:\smpq\mu]$, as $\la$ and $\mu$ vary.  The most basic results on this problem concern the dominance order.  If $\la$ and $\mu$ are partitions, we say that $\la$ \emph{dominates} $\mu$ (and write $\la\dom\mu$) if for each $i\gs1$ we have
\[\la_1+\dots+\la_i\gs\mu_1+\dots+\mu_i.\]
Now we have the following.

\begin{propn}\label{basicdecomp}\thmcite{dija}{Corollaries 4.12 \& 4.14}
Suppose $\la$ and $\mu$ are partitions of $n$, with $\mu$ $e$-regular.  Suppose $M$ is either the permutation module $\perq\la$ or the Specht module $\speq\la $.  If $[M:\smpq\mu]>0$, then $\mu\dom\la$.  If $\mu=\la$, then $[M:\smpq\mu]=1$.
\end{propn}

\subsection{Conjugate partitions and duality}\label{conjsec}

Let $T_1,\dots,T_{n-1}$ be the standard generators of $\hfq$.  Let $\sharp:\hfq\rightarrow\hfq$ be the involutory automorphism sending $T_i$ to $q-1-T_i$, and let $*:\hfq\rightarrow\hfq$ be the anti-automorphism sending $T_i$ to $T_i$.  Given a module $M$ for $\hfq$, define $M^\sharp$ to be the module with the same underlying vector space and with action
\[h\cdot m = h^\sharp m,\]
and define $M^{*}$ to be the module with underlying vector space dual to $M$ and with $\hfq$-action given by
\[h\cdot f(m) = f(h^{*}m).\]

We can describe the effect of these functors on Specht modules, using conjugate partitions.

\begin{propn}
Suppose $\la$ is a partition, and let $\la'$ denote the conjugate partition.  Then
\[(S^\la_{\bbf,q})^\sharp\equiv (S^{\la'}_{\bbf,q})^{*}.\]
\end{propn}

\begin{pf}
This is the result of \cite[Exercise 3.14(iii)]{mathbook}.  Although Mathas's definition of Specht modules is different from ours, his description of the relationship between the two definitions \cite[p.\ 38]{mathbook} ensures that the result holds for our Specht modules also.
\end{pf}

This has the following immediate corollary, which will be very useful.

\begin{cory}\label{conjirr}
Suppose $\la$ is a partition.  Then $S^\la_{\bbf,q}$ is irreducible if and only if $S^{\la'}_{\bbf,q}$ is.
\end{cory}

\subsection{Decomposition maps and adjustment matrices}\label{adjsec}

In this section we quote a result which will allow us to assume that our underlying field $\bbf$ is the field of rational numbers.  The theorem we shall state was proved by Geck in \cite{geck1}, and arises from a consideration of \emph{decomposition maps} between Iwahori--Hecke algebras defined over different rings.  (For an introduction to decomposition maps, see Geck's article \cite{geck2}.)  The theorem is most conveniently stated in terms of the \emph{decomposition matrix} of $\hfq$; this is the matrix with rows indexed by the partitions of $n$ and columns indexed by the $e$-regular partitions, and with the $(\la,\mu)$-entry being the decomposition number $[\speq\la :\smpq\mu]$.

\begin{thm}\label{adjmat}
Let $\zeta$ be a primitive $e$th root of unity in $\bbc$.  Let $D(\bbc,\zeta)$ and $D(\bbf,q)$ denote the decomposition matrices of $\hhh{\bbc,\zeta}(\sss n)$ and $\hfq$ respectively.  Then there is a square matrix $A$ with non-negative integer entries and with $1$s on the diagonal, such that
\[D(\bbf,q)=D(\bbc,\zeta)A.\]
\end{thm}

The important consequence of this theorem from our point of view is that for any field $\bbf$ and any partition $\la$, the Specht module $\speq\la$ has at least as many composition factors as $S^\la_{\bbc,\zeta}$; in particular, if $S^\la_{\bbc,\zeta}$ is reducible, then so is $\speq\la$.

\subsection{Cores and blocks}

Here we give the classification of the blocks of $\hsn$; this is based on the combinatorics of \emph{dominoes}.  Define a domino to be a pair of horizontally or vertically adjacent nodes in $\bbn^2$.  A \emph{removable domino} of a partition $\la$ is a domino contained in $[\la]$ such that the removal of this domino leaves the Young diagram of a partition.  The \emph{core} of $\la$ is the partition obtained by repeatedly removing removable dominoes until there are no more.  This partition has the form $(l,l-1,\dots,2,1)$ for some $l\gs0$, and is independent of the way is which the removable dominoes are chosen at each stage.  The \emph{weight} of $\la$ is the number of dominoes removed to obtain the core.

\begin{eg}
Let $\la=(6^2,5,2,1)$.  Then $\la$ has core $(3,2,1)$ and weight $7$, as we can see from the following diagram.
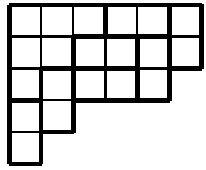
\[\begin{picture}(60,72)
\put(0,48){\line(1,0){24}}
\put(0,36){\line(1,0){12}}
\put(12,36){\line(0,1){24}}
\put(24,48){\line(0,1){12}}
\put(0,12){\line(1,0){12}}
\put(12,24){\line(1,0){12}}
\put(36,24){\line(0,1){24}}
\put(48,48){\line(0,1){12}}
\put(48,36){\line(1,0){12}}
\put(60,48){\line(1,0){12}}
\linethickness{0.5mm}\put(0,24){\line(0,1){36}}
\linethickness{0.5mm}\put(0,60){\line(1,0){36}}
\put(0,0){\vdom}
\put(12,12){\vdom}
\put(48,24){\vdom}
\put(60,36){\vdom}
\put(24,24){\hdom}
\put(24,36){\hdom}
\put(36,48){\hdom}
\end{picture}
\]
\end{eg}

Now we can address the blocks of $\hsn$.  Because the Specht modules are the cell modules arising from a particular cellular basis of $\hsn$, it follows that each Specht module lies entirely within one block of $\hsn$ \cite[Corollary 2.22]{mathbook}; we abuse notation by saying that a partition $\la$ lies in a block $B$ if $\spe\la$ lies in $B$.  Each block must contain at least one Specht module (because every simple module occurs as a composition factor of some Specht module), and so a classification of the blocks of $\hsn$ may be described by giving the appropriate partition of the set of partitions of $n$.  This is done as follows.

\begin{thm}\thmcite{mathbook}{Corollary 5.38}\label{blockclassn}
Suppose $\la$ and $\mu$ are two partitions of $n$.  Then $\la$ and $\mu$ lie in the same block of $\hsn$ if and only if $\la$ and $\mu$ have the same core.  Hence if $\mu$ is $2$-regular, then $[\spe\la :\smp\mu]=0$ unless $\la$ and $\mu$ have the same core.
\end{thm}

As a consequence of this theorem, we can define the weight and core of a block $B$, meaning the the weight and core of any partition labelling a Specht module lying in $B$.

\subsection{Rouquier blocks}


Suppose $B$ is a block of $\hsn$, with core $\nu=(l,l-1,\dots,1)$ and weight $w$.  We say that $B$ is \emph{Rouquier} if $w\ls l+1$.  Rouquier blocks are very useful because we have an explicit formula for their decomposition numbers in the case where the underlying field has infinite characteristic.  We now describe these results, following \cite{jm-1} (where a Rouquier block is referred to as a `block with an enormous $2$-core').

Suppose $B,\nu,w$ are as above, with $w\ls l+1$.  If $\la$ is a partition in $B$, then there is a unique way to partition the set $[\la]\setminus[\nu]$ into dominoes.  Given this partition, we define $\la^h_i$ to be the number of horizontal dominoes in row $i$ of $[\la]\setminus[\nu]$, and we define $\la^v_i$ to be the number of vertical dominoes in column $i$ of $[\la]\setminus[\nu]$, for $i\gs1$.  Then $\la^h$ and $\la^v$ are partitions, with $|\la^h|+|\la^v|=w$.  Moreover, $\la$ is uniquely specified by $\nu,\la^h,\la^v$, and in fact given any pair $\sigma,\tau$ of partitions with $|\sigma|+|\tau|=w$, there is a partition $\mu$ in $B$ with $\mu^h=\sigma$, $\mu^v=\tau$.

\begin{eg}
Suppose $\la=(13,8,7,4,3,2,1^5)$.  Then $\la$ has weight $7$ and core $(7,6,5,4,3,2,1)$, so lies in a Rouquier block.  The Young diagram $[\la]$ may be drawn as follows, and we see that $(\la^h,\la^v)=\left((3,1^2),(2)\right)$:
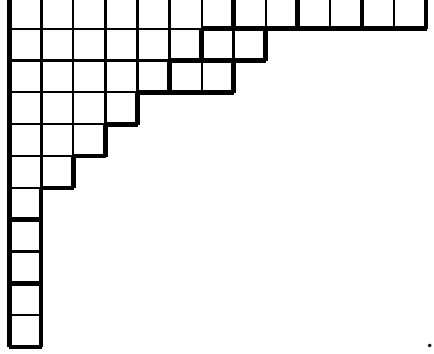
\[\begin{picture}(156,132)
\multiput(0,12)(0,24)3{\line(1,0){12}}
\put(0,72){\line(1,0){24}}
\put(0,84){\line(1,0){36}}
\put(0,96){\line(1,0){48}}
\put(0,108){\line(1,0){60}}
\put(0,120){\line(1,0){72}}
\put(12,60){\line(0,1){72}}
\put(24,72){\line(0,1){60}}
\put(36,84){\line(0,1){48}}
\put(48,96){\line(0,1){36}}
\put(60,108){\line(0,1){24}}
\put(72,96){\line(0,1){12}}
\put(84,108){\line(0,1){12}}
\multiput(72,120)(24,0)4{\line(0,1){12}}
\put(0,0){\vdom}
\put(0,24){\vdom}
\multiput(84,120)(24,0){3}{\hdom}
\put(72,108){\hdom}
\put(60,96){\hdom}
\multiput(12,48)(12,12){4}{\line(0,1){12}}
\multiput(12,60)(12,12){4}{\line(1,0){12}}
\linethickness{0.5mm}
\put(0,48){\line(0,1){84}}
\put(0,132){\line(1,0){84}}
\multiput(12,48)(12,12)4{\put(0,0){\line(0,1){12}}\put(0,12){\line(1,0){12}}}
\end{picture}.\]

\end{eg}

So we may label a partition $\la$ in a Rouquier block $B$ by the pair $(\la^h,\la^v)$; we remark that the pair $(\la^h,(\la^v)')$ is often referred to as the \emph{$2$-quotient} of $\la$.  It is easy to see that $\la$ is $2$-regular if and only if $\la^v=\varnothing$, while $\la$ is $2$-restricted if and only if $\la^h=\varnothing$.  Now we can describe the decomposition numbers for a Rouquier block in infinite characteristic; given any partitions $\alpha,\beta,\gamma$ with $|\alpha|=|\beta|+|\gamma|$, let $c^\alpha_{\beta\gamma}$ be the corresponding Littlewood--Richardson coefficient (see \cite{fult}).

\begin{thm}\thmcite{jm-1}{Theorem 2.5}\label{rouqdec}
Suppose $\nchar(\bbf)=\infty$ and $B$ is a Rouquier block of $\hsn$.  Suppose $\la$ and $\mu$ are partitions in $B$, with $\mu$ $2$-regular, and let $(\la^h,\la^v)$ and $(\mu^h,\varnothing)$ be the corresponding pairs of partitions.  Then
\[[\spe\la :\smp\mu]=c^{\mu^h}_{\la^h\la^v}.\]
\end{thm}

The only corollary we need from this is the following.

\begin{cory}\label{rouqcory}
Suppose $B$ is a Rouquier block of $\hsn$, and $\la$ is a partition in $B$ which is neither $2$-regular nor $2$-restricted.  Then $\spe\la $ is reducible.
\end{cory}

\begin{pf}
By the results of \S\ref{adjsec}, we may assume that $\nchar(\bbf)=\infty$.  Since $\la$ is neither $2$-regular nor $2$-restricted, the labelling partitions $\la^h$ and $\la^v$ are both non-empty.  It is well known and easy to prove from the definition of Littlewood--Richardson coefficients that if $\beta$ and $\gamma$ are non-empty partitions, then there are at least two partitions $\alpha$ for which $c^\alpha_{\beta\gamma}>0$, and now the result follows from Theorem \ref{rouqdec}.
\end{pf}

\subsection{Alternating partitions}

The question of irreducibility of Specht modules labelled by $2$-regular partitions has been settled for some time.  We shall need this result later, so we quote it here, concentrating for simplicity on the case where $\nchar(\bbf)=\infty$.  Say that a partition is \emph{alternating} if for every $i$ either $\la_i+\la_{i+1}$ is odd or $\la_{i+1}=0$.

\begin{thm}\thmcite{mathbook}{Theorem 5.42} \label{Carter}
Suppose $\nchar(\bbf)=\infty$, and $\la$ is a $2$-regular partition.  Then $\spe\la $ is irreducible if and only if $\la$ is alternating.
\end{thm}

We combine this with the following simple lemma.

\begin{lemma}\label{isalternating}
Suppose $\la$ is a partition with core $(l-1,l-2,\dots,1)$ for some $l\gs1$, and that $\la_{l+1}=0$.  Then $\la$ is alternating.
\end{lemma}

\begin{pf}
We use induction to prove a stronger statement, namely that $\la_i\equiv i+l \pmod 2$ for $1\ls i\ls l$.  The induction is on the weight $w$ of $\la$; if $w=0$ then $\la$ is the core $(l-1,l-2,\dots,1)$, which certainly has the required property.  For the inductive step, suppose that $w>0$ and let $\nu$ be a partition obtained by removing a domino from $[\la]$.  Then $\nu$ satisfies the hypotheses of the lemma, so by induction we have $\nu_i\equiv i+l\pmod 2$ for $i=1,\dots,l$.  In particular, we have $\nu_i\neq\nu_{i+1}$ if $1\ls i\ls l-1$.  This means that the domino added to $[\nu]$ to obtain $[\la]$ must be horizontal; for if it were vertical, consisting of the nodes $(i,j)$ and $(i+1,j)$ say, then we would have $\nu_i=\nu_{i+1}(=j-1)$, and hence $i\gs l$.  But this would give $\la_{l+1}>0$, contradicting our assumptions.  So the added domino is horizontal, and hence $\la_i\equiv\nu_i\pmod 2$ for all $i$, which gives the required conclusion.
\end{pf}

When we combine Lemma \ref{isalternating} with Proposition \ref{basicdecomp}, the following result is immediate.

\begin{cory}\label{altcory}
Suppose $\nchar(\bbf)=\infty$, and $\mu$ is a $2$-regular partition with core $(l-1,l-2,\dots,1)$.  If $\la$ is a partition such that $[\spe\la :\smp\mu]>0$ and $\la_{l+1}=0$, then $\la=\mu$.
\end{cory}

\subsection{Ladders and James's regularisation theorem}\label{ladsec}

We saw above that when $\la$ is $2$-regular, the simple module $\smp\la$ occurs as a composition factor of $\spe\la $ with multiplicity $1$.  James has given an extension of this result to the case where $\la$ is not $2$-regular, giving an explicit simple module which occurs exactly once as a composition factor of $\spe\la $.  This was done in \cite{j1} for symmetric groups, and extended to the general case in \cite{j10}.  This is very useful from the point of view of classifying irreducible Specht modules, since if $\spe\la $ is irreducible, the theorem tells us which irreducible module $D^\mu_{\bbf,-\negthinspace1}$ is isomorphic to $\spe\la $.

Recall the definition of ladders from \S\ref{irrsumm}.  Given a partition $\la$, it is easily seen that $\la$ is $2$-regular if and only if for each $l$ the nodes in the $l$th ladder of $\la$ are as high as possible, i.e.\ $\call_l(\la) = \left\{(1,l),(2,l-1),\dots,(s,l+1-s)\right\}$ for some $s$.  Furthermore, for any $\la$ we may obtain a $2$-regular partition by moving the nodes in each ladder of $\la$ to the topmost positions in that ladder; this $2$-regular partition is called the \emph{regularisation} of $\la$, written $\la\reg$.

\begin{eg}
Let $\la=(4,2^3)$.  Then we have
\begin{align*}
\call_l(\la) &= \call_l \text{ for $l\ls 3$},\\
\call_4(\la) &= \{(1,4),(3,2),(4,1)\},\\
\call_5(\la) &= \{(4,2)\},\\
\call_l(\la) &=\emptyset\text{ for $l\gs6$}.
\end{align*}
By replacing $\call_4(\la)$ with $\{(1,4),(2,3),(3,2)\}$ and $\call_5(\la)$ with $\{(1,5)\}$, we obtain $\la\reg=(5,3,2)$:
\[\la = \raisebox{-24pt}{\begin{picture}(48,48)\multiput(0,36)(12,0)4{\bx}\multiput(0,0)(12,0)2{\multiput(0,0)(0,12)3{\bx}}\end{picture}};\qquad\la\reg=\raisebox{-12pt}{\begin{picture}(60,36)\multiput(0,24)(12,0)5\bx\multiput(0,12)(12,0)3\bx\multiput(0,0)(12,0)2\bx\end{picture}}.\]
\end{eg}

\begin{thm}\thmcite{j10}{Theorem 6.21}\label{regthm}
Suppose $\la,\mu$ are partitions of $n$, with $\mu$ $2$-regular.  If $[\spe\la :\smp\mu]>0$, then $\mu\dom\la\reg$.  Furthermore, $[\spe\la :\smp{\la\reg}]=1$.
\end{thm}

We note the following corollary concerning homomorphisms.

\begin{cory}\label{homcory}
Suppose $\la,\mu$ are partitions of $n$.  Suppose that either:
\begin{enumerate}
\item
$\la\reg\ndom\mu\reg$, and there exists a non-zero homomorphism from $\spe\mu$ to $\spe\la$; or
\item
$\la\reg\ndom\mu$, and there exists a non-zero homomorphism from $\per\mu$ to $\spe\la $.
\end{enumerate}
Then $\spe\la $ is reducible.
\end{cory}

\begin{pf}
This follows immediately from Theorem \ref{regthm} and Proposition \ref{basicdecomp}.
\end{pf}



%

\subsection{A useful lemma}\label{213sec}

In this section we recall a very useful result; the general version of this result was the main tool used in \cite{slred} for proving the reducibility of Specht modules.

Given a partition $\la$ and $i\in\{0,1\}$, define the partition $\la^{-i}$ by removing all removable nodes of residue $i$.  Then we have the following, which is a corollary of \cite[Lemma 2.13]{bk}.

\begin{lemma}\label{213}
Suppose $\la$ is a partition and $i\in\{0,1\}$.  Then $\spe\la$ has at least as many composition factors as $\spe{\la^{-i}}$.  In particular, if $\spe{\la^{-i}}$ is reducible, then so is $\spe\la$.
\end{lemma}

We also need a `dual' result to this.  Given $\la$ and $i\in\{0,1\}$, define $\la^{+i}$ by adding all addable nodes of residue $i$.  Then the following result may be proved in exactly the same way as Lemma \ref{213}.

\begin{lemma}\label{213dual}
Suppose $\la$ is a partition and $i\in\{0,1\}$.  Then $\spe\la$ has at least as many composition factors as $\spe{\la^{+i}}$.  In particular, if $\spe{\la^{+i}}$ is reducible, then so is $\spe\la$.
\end{lemma}

To use the latter result in inductive proofs, it will be helpful to have the following lemma, in which we write $w(\la)$ for the weight of a partition $\la$.

\begin{lemma}\label{wtdown}
Suppose $\la$ is a partition and $i\in\{0,1\}$.  Then $w(\la^{+i})\ls w(\la)$.
\end{lemma}

\begin{pf}
A much more general result is proved in \cite[Lemma 3.6]{mfweight}.  
\end{pf}

\sect{Homomorphisms}\label{homsec}

In this section we recall some results concerning the existence of homomorphisms between permutation modules and Specht modules, and we use these to construct homomorphisms in particular cases; we also recall the second author's analogue of a special case of the Carter--Payne theorem for homomorphisms between Specht modules.  In conjunction with Corollary \ref{homcory}, these results will be useful for proving reducibility of Specht modules.

\subsection{Homomorphisms from permutation modules to Specht modules}

Suppose $\mu$ is a composition of $n$.  A \emph{$\mu$-tableau} is a function $T$ from $[\mu]$ to $\bbz_{\gs0}$; given a tableau $T$, we write $T(i,j)$ instead of $T((i,j))$, and we usually depict $T$ by drawing the Young diagram $[\mu]$, and filling each node with with its image under $T$.  Given $i,j\gs1$, we write $T_{i,j}$ for the number of entries equal to $j$ in row $i$ of $T$.  If $\la$ is another composition of $n$, then we say that a $\mu$-tableau $T$ has \emph{content} $\la$ if for every $i$ there are exactly $\la_i$ nodes mapped to $i$ by $T$.

\begin{eg}
Let $\mu=(5,3,1)$ and $\la=(4,3,2)$.  Then the tableau
\[T=\raisebox{-15pt}{\young(11133,122,2)}\]
is a $\mu$-tableau of content $\la$.  The values $T_{i,j}$ are given by the following matrix:
\[\begin{array}{c|ccc}
&1&2&3\\\hline
1&3&0&2\\2&1&2&0\\3&0&1&0\end{array}.\]
\end{eg}

For each $\mu$-tableau $T$ of content $\la$, Dipper and James define a homomorphism $\Theta_T:M^\mu_{\bbf,q}\rightarrow M^\la_{\bbf,q}$, for any $\bbf,q$.  The homomorphisms $\Theta_T$ and $\Theta_U$ are equal if $T$ and $U$ are \emph{row equivalent}; that is, $U_{i,j}=T_{i,j}$ for each $i,j$.  We say that $T$ is \emph{row standard} if the entries in $T$ are weakly increasing along rows, and we write $\calt(\mu,\la)$ for the set of row standard $\mu$-tableaux of content $\la$.  Dipper and James prove that the set
\[\big\{\Theta_T\ \big|\ T\in \calt(\mu,\la)\big\}\]
is a basis for the space of homomorphisms from $M^\mu_{\bbf,q}$ to $M^\la_{\bbf,q}$.

A particular set of homomorphisms $\Theta_T$ can be used to give a convenient characterisation of the Specht module.  Suppose $\la$ is a partition, and suppose $d,t\gs1$ are such that $t\ls\la_{d+1}$.  Define the composition $\la^{d,t}$ by
\[\la^{d,t}_i = \begin{cases}
\la_d+t & (i=d)\\
\la_{d+1}-t & (i=d+1)\\
\la_i & (\text{otherwise}).
\end{cases}\]
Then there is a unique tableau $A\in\calt(\la,\la^{d,t})$ with the property that $A(i,j)=i$ for all $(i,j)\in[\la]$ with $i\neq d+1$.  We write $\psi^{d,t}$ for the homomorphism $\Theta_A:M^\la_{\bbf,q}\rightarrow M^{\la^{d,t}}_{\bbf,q}$.  (Warning: in \cite{slcp} and elsewhere, the map $\psi^{d,t}$ is written as $\psi_{d,\la_{d+1}-t}$.  Our convention is more convenient here.)

Now the Specht module can be characterised as follows.

\begin{thm}\thmcite{dija}{Theorem 7.5}\label{kerint}
Suppose $\la$ is a partition.  Then
\[S^\la_{\bbf,q}=\bigcap_{d\gs1}\bigcap_{1\ls t\ls\la_{d+1}}\ker\psi^{d,t}.\]
\end{thm}

This theorem is known as the kernel intersection theorem, and is very useful for constructing and classifying homomorphisms $M\rightarrow S^\la_{\bbf,q}$, when $M$ is a module for which one knows all homomorphisms $M\rightarrow M^\la_{\bbf,q}$; specifically, the homomorphisms $M\rightarrow S^\la_{\bbf,q}$ are precisely the homomorphisms $\Theta:M\rightarrow M^\la_{\bbf,q}$ such that $\psi^{d,t}\circ\Theta=0$ for all $d,t$.  This approach has been particularly exploited by the second author, using an explicit description of the composition $\psi^{d,t}\circ\Theta_T$, when $T\in\calt(\mu,\la)$.  Before we can give this result, we need to give a very brief account of quantum binomial coefficients.  For any non-negative integer $a$, we define the \emph{quantum integer} $[a]=1+q+q^2+\dots+q^{a-1}$, and the \emph{quantum factorial} $[a]! = [1][2]\dots[a]$.  Now for $0\ls b\ls a$ the \emph{quantum binomial coefficient} is defined to be
\[\dqbinom ab = \frac{[a]!}{[b]![a-b]!}.\]
Of course, if $q=1$ then this coincides with the usual binomial coefficient $\binom ab$.

The only property of quantum binomial coefficients we need is the following, which is well-known.

\begin{lemma}\label{qbc}
Suppose $q=\mo$, and $0\ls b\ls a$.  Then
\[\dqbinom ab = \begin{cases}
\dbinom{\lfloor a/2\rfloor}{\lfloor b/2\rfloor} & (\text{if $a$ is odd or $b$ is even})\\
\ \ \ \ \ 0 & (\text{if $a$ is even and $b$ is odd}).
\end{cases}\]
\end{lemma}

Suppose $\mu$ and $\la$ are two partitions, and $d,t$ are chosen as above.  Given $T\in\calt(\mu,\la)$, let $\calv_T\subseteq\calt(\mu,\la^{d,t})$ be the set of row-standard tableaux $V$ with the property that for each $(i,j)\in[\mu]$ either $V(i,j)=T(i,j)$ or $V(i,j)=d=T(i,j)-1$.  (In other words, $V$ is a row-standard tableau obtained from $T$ by replacing $t$ of the entries equal to $d+1$ with $d$s.)

Given $V\in\calv_T$, define
\[x = \sum_{i\gs1}\left((V_{i,d}-T_{i,d})\sum_{k>i}T_{k,d}\right),\]
and set
\[b^{(q)}_{TV} = q^x\prod_{i\gs1}\dqbinom{V_{i,d}}{T_{i,d}},\]
considered as an element of the field $\bbf$.  Now we have the following statement.

\begin{propn}\thmcite{slcp}{Proposition 2.14}\label{psicomp}
Suppose $\la,\mu,T,d,t$ are as above.  Then
\[\psi^{d,t}\circ\Theta_T = \sum_{V\in\calv_T}b^{(q)}_{TV}\Theta_V.\]
\end{propn}

We shall use this to prove the following proposition.

\begin{propn}\label{mattshom}
Suppose $\nu$ and $\xi$ are partitions; put $l=\nu'_1$, and suppose that $\xi_{l-1}\gs l$.  Define partitions $\la,\mu$ by
\begin{align*}
\la_i &= \xi_i+2\nu_i,\\
\mu_i &= \xi'_i+2\nu_i
\end{align*}
for all $i\gs1$.  Then there is a non-zero $\hh$-homomorphism from $\per\mu$ to $\spe\la$.
\end{propn}

\begin{eg}
Put $\xi=(2^4)$ and $\nu=(1^2)$.  Then we have $l=2$, so that $\xi_{l-1}\gs l$, and so when $q=\mo$ there is a non-zero homomorphism from $\per\mu$ to $\spe\la$, where $\la=(4^2,2^2)$ and $\mu=(6^2)$.  We have $\la\reg=(5,4,2,1)$, which does not dominate $\mu$, and so by Corollary \ref{homcory} we deduce that $\spe{(4^2,2^2)}$ is reducible.
\end{eg}

\begin{pfof}{Proposition \ref{mattshom}}
Using the kernel intersection theorem, we need to construct a linear combination
\[\theta = \sum_{T\in\calt(\mu,\la)}c_T\Theta_T\]
such that $\psi^{d,t}\circ\theta=0$ for all $d,t$.  To do this, we begin with the semi-simple Iwahori--Hecke algebra $\hhh\bbq=\hhh{\bbq,1}(\sss m)$, i.e.\ the group algebra $\bbq\sss m$, where $m=|\nu|$.  The module $M^\nu_{\bbq,1}$ contains the Specht module $S^\nu_{\bbq,1}$ as a submodule; since $\hhh\bbq$ is semi-simple, $S^\nu_{\bbq,1}$ is also a quotient of $M^\nu_{\bbq,1}$, so there is a non-zero $\hhh\bbq$-homomorphism $\phi:M^\nu_{\bbq,1}\rightarrow S^\nu_{\bbq,1}$.  Regarding this as a homomorphism from $M^\nu_{\bbq,1}$ to itself, and using the fact that the homomorphisms $\Theta_T$ for $T\in\calt(\nu,\nu)$ span the space of such homomorphisms, we may write
\[\phi = \sum_{T\in\calt(\nu,\nu)}a_T\Theta_T\]
with $a_T\in\bbq$.  By re-scaling, we may assume that the $a_T$ are coprime integers.

The fact that the image of $\phi$ lies in the Specht module $S^\nu_{\bbq,1}$ implies that
\[\sum_{T\in\calt(\nu,\nu)}a_T\psi^{d,\tau}\circ\Theta_T=0\]
for all $d,\tau$ with $1\ls\tau\ls\nu_{d+1}$.  By Proposition \ref{psicomp}, this means that for each pair $d,\tau$ we have
\[\sum_{T\in\calt(\nu,\nu)}a_T\left(\sum_{V\in\calv_T}b^{(1)}_{TV}\Theta_V\right)=0.\]
Since the set $\left\{\Theta_V\ \left|\ V\in\calt(\nu,\nu^{d,\tau})\right.\right\}$ is linearly independent, this says that for each $V\in\calt(\nu,\nu^{d,\tau})$, the sum
\[\sum_{T\in\calt(\nu,\nu)\ \mid\  V\in\calv(T)}a_Tb^{(1)}_{TV}\]
vanishes.

Now we construct a homomorphism $\theta:\per\mu\rightarrow \per\la$, whose image we claim lies in $\spe\la $.  For each $T\in\calu(\nu,\nu)$, let $\hat T$ be the row-standard $\mu$-tableau given by
\[\hat T_{i,j} = \begin{cases}
2T_{i,j}+1 & (j\ls \xi'_i)\\
2T_{i,j} & (j>\xi'_i)
\end{cases}\]
for each $i,j$.
  Then we have $\hat T\in\calt(\mu,\la)$, and we define
\[\theta = \sum_{T\in\calt(\nu,\nu)} a_T\Theta_{\hat T}.\]
(We are committing a minor abuse of notation here: by $a_T$, we really mean the image of $a_T$ in $\bbf$, which is well-defined since we are assuming that each $a_T$ is an integer.)  The tableau $T$ is easily recovered from $\hat T$, so the tableaux $\hat T$ are distinct, and therefore the homomorphisms $\Theta_{\hat T}$ are linearly independent; since we assume that the integers $a_T$ are coprime, this implies that $\theta$ is non-zero.  (Alternatively, one can use the results of \S\ref{adjsec} and assume throughout that $\bbf=\bbq$.)

Fix a pair $d,t$ with $1\ls t\ls \la_{d+1}$.  Using the kernel intersection theorem and Proposition \ref{psicomp}, our task is to show that the sum
\[\sum_{T\in\calt(\nu,\nu)}a_T\sum_{V\in\calv_{\hat T}}b^{(\mo)}_{\hat TV}\Theta_V\]
equals zero.

\clam
Suppose $T\in\calt(\nu,\nu)$ and $V\in\calv_{\hat T}$, and define
\[\beta_i=\left|\left\{j\ \left|\ V(i,j)\neq \hat T(i,j)\right.\right\}\right|\]
as above.  If $\beta_i$ is odd for any $i$, then $b^{(\mo)}_{\hat TV}=0$.
\prof
Fix $i$ such that $\beta_i$ is odd.  Then we claim that $\hat T_{i,d}$ is odd.  This will then imply that the integer $y_i=V_{i,d}$ is even, so that $\qbinom{y_i}{\beta_i}$ equals zero, and hence $b^{(\mo)}_{\hat TV}$ is zero.

The fact that $\beta_i>0$ means that $\hat T_{i,d+1}>0$.  Write $l=\nu'_1$ as above, and suppose first that $i>l$.  Then by construction the entries in row $i$ of $\hat T$ are $1,2,\dots,\xi'_i$ each occurring once; since $d+1$ occurs, we have $d+1\ls \xi'_i$, so that $\hat T_{i,d}=1$.

Alternatively, suppose that $i\ls l$.  Then we claim that $\xi'_i\gs d$, which will mean that
\[\hat T_{i,d} = 1+2T_{i,d},\]
which is odd.  We are given that $\xi_{l-1}\gs l$, i.e.\ $\xi'_l\gs l-1$, and hence $\xi'_i\gs l-1$.  So if $d\ls l-1$ we are done.  If $d>l-1$ then $d+1>l$, so $T_{i,d+1}=0$; but by assumption there is an entry equal to $d+1$ in row $i$ of $\hat T$, and hence we must have $\xi'_i\gs d+1$.
\malc

As a consequence of the claim, we need only consider those pairs $(T\in\calt(\nu,\nu),V\in\calv_{\hat T})$ for which $\beta_i$ is even for each $i$.  Given such a pair, this condition implies that $t=\sum_i\beta_i$ is even and there is a unique $W\in\calt(\nu,\nu^{d,t/2})$ such that $V=\hat W$ (where we define $\hat W$ analogously to $\hat T$).  Furthermore, for each such $V$ and each $T\in\calt(\nu,\nu)$ we have $V\in\calv_{\hat T}$ if and only if $W\in\calv_T$, and if this is the case then it is easy to calculate that $b^{(\mo)}_{\hat TV}=b^{(1)}_{TV}$, using Lemma \ref{qbc}.  Now the result follows.
\end{pfof}

\begin{eg}
Let $\xi,\nu$ be as in the previous example.  The two tableaux in $\calt(\nu,\nu)$ are
\setlength{\raiser}{9pt}
\[T_1 = \raisebox{-\raiser}{\young(1,2)},\qquad T_2=\raisebox{-\raiser}{\young(2,1)},\]
and a non-zero homomorphism from $M^\mu_{\bbq,1}$ to $S^\mu_{\bbq,1}$ is given by
\[\phi = \Theta_{T_1}-\Theta_{T_2}.\]
So a homomorphism from $\per{(6^2)}$ to $\spe{(4^2,2^2)}$ is given by
\[\theta = \Theta_{\hat T_1}-\Theta_{\hat T_2},\]
where
\[\hat T_1 = \raisebox{-\raiser}{\young(111234,122234)},\qquad\hat T_2=\raisebox{-\raiser}{\young(122234,111234)}.\]
\end{eg}

\subsection{Homomorphisms between Specht modules}

We now quote a result due to the second author which gives the existence of non-zero homomorphisms between Specht modules under certain circumstances.  This is a generalisation to Iwahori--Hecke algebras of a special case of the Carter--Payne Theorem \cite{cp}.  Recall the notion of the residue of a node from \S\ref{213sec}.

\begin{thm}\thmcite{slcp}{Theorem 4.1.1}\label{slcp}
Suppose $\la$ is a partition, and that $\la$ has an addable node $(i,\la_i+1)$ and a removable node $(j,\la_j)$ of the same residue, with $i<j$.  Let $\mu$ be the partition obtained by adding the node $(i,\la_i+1)$ and removing the node $(j,\la_j)$.  Then there exists a non-zero homomorphism from $\spe\mu$ to $\spe\la $.
\end{thm}

This result will be helpful in conjunction with Corollary \ref{homcory}.  A particular application is the following.

\begin{propn}\label{ladhom}
Suppose $\la$ is a partition.  Suppose $\la$ has
\begin{itemize}
\item
an addable node $(i,\la_i+1)$ lying in ladder $\call_m$, and
\item
a removable node $(j,\la_j)$ lying in ladder $\call_l$,
\end{itemize}
where $m>l$ and $l\equiv m\pmod 2$.  Then $\spe\la $ is reducible.
\end{propn}

\begin{pf}
By replacing $\la$ with its conjugate if necessary and appealing to Corollary \ref{conjirr}, we may assume that $i<j$.  Since $l\equiv m\pmod2$ and the nodes in ladder $\call_l$ all have residue $(l+1)\pmod 2$, the addable node $(i,\la_i+1)$ and the removable node $(j,\la_j)$ both have the same residue.  So if we define $\mu$ as in Theorem \ref{slcp}, then there is a non-zero homomorphism from $\spe\mu$ to $\spe\la $.  By Corollary \ref{homcory}, it suffices to show that $\la\reg\ndom\mu\reg$; this follows from \cite[Lemma 2.1]{mfreduc}, given the assumption that $m>l$.
\end{pf}

\sect{The Fock space and canonical bases}\label{focksec}

Now we introduce the Fock space, which is our most powerful tool.  In fact, via Ariki's Theorem, this theory provides an algorithm for computing the decomposition matrix of $\hsn$ completely when $\nchar(\bbf)=\infty$.  However, it does not seem easy to use this algorithm to decide the reducibility of Specht modules, and our application of the Fock space will be less direct.

Let $v$ be an indeterminate over $\bbq$, and let $\fkh=\bbq h_0\oplus\bbq h_1\oplus\bbq D$ be a three-dimensional vector space.  In this section we work with the quantum group $\calu=\calu_v(\widehat{\mathfrak{sl}}_2)$, which may be realised as the associative algebra over $\bbq$ with generators $e_0,e_1,f_0,f_1$ and $v^h\ (h\in\fkh)$, subject to well-known relations; these may be found in \cite{llt}, which is an excellent background reference for this section.

Define the \emph{Fock space} to be the $\bbq(v)$-vector space $\calf$ with a basis $\{s(\la)\}$ indexed by the set of all partitions.  Let $\lan\, ,\,\ran$ be the inner product on $\calf$ for which the basis $\{s(\la)\}$ is orthonormal.   The Fock space has the structure of a $\calu$-module, and has important connections to the representation theory of Iwahori--Hecke algebras.  It will suffice for our purposes to describe the action on $\calf$ of the `negative' generators $f_0,f_1$ and their `quantum divided powers'
\[f^{(a)}_i = \frac{f^a_i}{v^{a-1}+v^{a-3}+\dots+v^{3-a}+v^{1-a}}.\]
Fix $i\in\{0,1\}$ and $a\gs1$, and suppose $\la$ and $\mu$ are partitions.  Write $\la\stackrel{a:i}\lra\mu$ if the Young diagram for $\mu$ may be obtained from the Young diagram for $\la$ by adding $a$ addable nodes of residue $i$.  If this is the case, then for each $j$ such that $\la_j=\mu_j$ define
\[\epsilon_j = \begin{cases}
+\negthinspace1 & (\rez(j,\la_j+1)=i)\\
\mo & (\rez(j,\la_j+1)=1-i),
\end{cases}\]
and set
\[N(\la,\mu) = \sum_{j\ \mid\ \la_j=\mu_j}\epsilon_j\times\left(\text{no.\ of nodes of $[\mu]\setminus[\la]$ below row $j$}\right).\]

Now the action on $\calf$ of the quantum divided power $f_i^{(a)}$ is given by $\bbq(v)$-linear extension of
\[f^{(a)}_is(\la) = \sum_{\mu\ \mid\ \la\stackrel{a:i}\lra\mu}v^{N(\la,\mu)}s(\mu).\]

Of particular importance is the submodule of $\calf$ generated by the vector $s(\varnothing)$.  This submodule is isomorphic to (and therefore often identified with) the irreducible highest-weight representation $M(\Lambda_0)$ of $U$.  This representation is equipped with a $\bbq(v+v^{\mo})$-linear map called the \emph{bar involution}, which can be specified by the conditions
\[\ol{s(\varnothing)}=s(\varnothing)\]
and
\[\ol{f_i(m)} = f_i(\ol m)\]
for $i\in\{0,1\}$ and $m\in M(\Lambda_0)$.

The bar involution allows us to define the \emph{canonical basis} of $M(\Lambda_0)$, via the following theorem.

\begin{thm}\label{canbas}
For each $2$-regular partition $\mu$ there is a unique element $G(\mu)$ of $M(\Lambda_0)$ with the properties
\begin{itemize}
\item
$\ol{G(\mu)}=G(\mu)$, and
\item
$G(\mu) = \sum_{\la}d_{\la\mu}(v)s(\la)$, where $d_{\la\mu}(v)$ is a polynomial in $v$, with $d_{\mu\mu}(v)=1$, and $d_{\la\mu}(v)$ divisible by $v$ for $\la\neq\mu$.
\end{itemize}
The set
\[\left\{G(\mu)\ \left|\ \mu\text{ a $2$-regular partition}\right.\right\}\]
is a $\bbq(v)$-basis of $M(\Lambda_0)$.
\end{thm}

Now we can state (a special case of) Ariki's Theorem, which gives the connection to the representation theory of Iwahori--Hecke algebras.

\begin{thm}\thmcite{ari}{Theorem 4.4}\label{arikithm}
Suppose $\la$ and $\mu$ are partitions of $n$, with $\mu$ $2$-regular, and let $d_{\la\mu}(v)=\lan G(\mu),s(\la)\ran$ as in Theorem \ref{canbas}.  Then
\[[S^\la_{\bbq,\mo}:D^\mu_{\bbq,\mo}]=d_{\la\mu}(1).\]
\end{thm}

In view of this theorem, the polynomials $d_{\la\mu}(v)$ are known as `$v$-decomposition numbers'.  It is known \cite[Theorem 6.28]{mathbook} that $d_{\la\mu}(v)$ has non-negative integer coefficients, and is zero unless $\la$ and $\mu$ have the same core and weight (and therefore the same size).  The non-negativity of the coefficients has the following obvious consequence, in conjunction with Theorems \ref{regthm} and \ref{arikithm}.

\begin{lemma}\label{mono}
Suppose $\la$ and $\mu$ are partitions of $n$, with $\mu$ $2$-regular.
\begin{enumerate}
\item
If $[\spe\la:\smp\mu]=0$, then $d_{\la\mu}(v)=0$.
\item
If $[\spe\la:\smp\mu]=1$ (in particular, if $\mu=\la\reg$), then $d_{\la\mu}(v)=v^s$ for some $s$.
\end{enumerate}
\end{lemma}

\begin{rmk}
In the case where $\mu=\la\reg$, the integer $s$ in Lemma \ref{mono} has been computed explicitly by the first author in \cite{mfreg}; however, we shall not need this result in the present paper.
\end{rmk}

Next we prove a crucial lemma, which enables us to use Fock space computations to prove reducibility of Specht modules.

\begin{lemma}\label{mctsl}
Suppose $\la$ is a partition, and suppose $X$ and $Y$ are bar-invariant elements of $M(\Lambda_0)$, such that
\[\lan X,s(\la)\ran=v^x,\qquad\lan Y,s(\la)\ran=v^y\]
for some $x\neq y$.  Then the Specht module $S^\la_{\bbq,\mo}$ is reducible.
\end{lemma}

\begin{pf}
Suppose not, and write $\nu=\la\reg$.  Then by Theorem \ref{regthm} and Lemma \ref{mono} we have $\lan G(\mu),s(\la)\ran=0$ for any $\mu\neq\nu$.  We write $X$ and $Y$ as linear combinations of canonical basis vectors
\[X = \sum_\mu\alpha_\mu(v) G(\mu),\qquad Y = \sum_\mu\beta_\mu(v) G(\mu);\]
since $X$ and $Y$ are bar-invariant, $\alpha_\mu(v)$ and $\beta_\mu(v)$ lie in $\bbq(v+v^{\mo})$ for each $\mu$.
Taking inner products with $s(\la)$ yields
\[\alpha_{\nu}(v)d_{\la\nu}(v)=v^x,\qquad \beta_{\nu}(v)d_{\la\nu}(v)=v^y.\]
This gives
\[v^y\alpha_{\nu}(v) = v^x\beta_{\nu}(v)\]
with $\alpha_{\nu}(v)$ and $\beta_{\nu}(v)$ non-zero, but this is impossible if $x\neq y$ and $\alpha_{\nu}(v)$ and $\beta_{\nu}(v)$ lie in $\bbq(v+v^{\mo})$.
\end{pf}

\begin{rmk}
In fact, Lemma \ref{mctsl} shows something rather stronger than the reducibility of $\spe\la$.  Let us say that $\la$ is \emph{homogeneous} if there is some $x$ such that every $v$-decomposition number $d_{\la\mu}(v)$ is either zero or a monomial of degree $x$.  According to popular conjectures relating $v$-decomposition numbers to the Jantzen filtration of the Specht module, the homogeneity of $\la$ ought to imply that the Specht module $S^\la_{\bbq,\mo}$ is completely reducible.

A weaker condition we might impose on $\la$ is that it is \emph{quasi-homogeneous}, meaning that there is some $x$ such that every $d_{\la\mu}(v)$ lies in $v^x.\bbq(v+v^{\mo})$.  As a representation-theoretic interpretation, we would speculate that quasi-homogeneity corresponds to the Specht module $\spe\la$ being self-dual.

What we have actually shown is that if the hypotheses of Lemma \ref{mctsl} are satisfied, then $\la$ is not homogeneous or even quasi-homogeneous.  It would be very interesting to classify homogeneous and quasi-homogeneous partitions, and the authors hope to be able to say something more about this in the future.
\end{rmk}

Our next step is to prove a result in which we phrase the action of a certain composition of powers of $f_0,f_1$ in a convenient form, in certain special cases.  Fix partitions $\mu$ and $\la$ with $[\mu]\subseteq[\la]$, and fix $x\in\{0,1\}$.  Suppose that the following condition holds:
\[\text{for each $i$ with $\la_i>0$, the node $(i,\mu_i+1)$ has residue $x$.}\tag*{($\ast$)}\]
In other words, $\mu_i\equiv i+x\pmod2$ whenever $\la_i>0$.  Note that this implies in particular that $\la'_1-\mu'_1\ls1$.

Define a sequence of partitions $\mu=\mu^0,\mu^1,\mu^2,\dots$ as follows: for $j\gs0$, $\mu^{j+1}$ is obtained from $\mu^j$ by adding all addable nodes that are contained in $[\la]$.  
We define $a_j=|\mu^j|-|\mu^{j-1}|$ for $j\gs 1$, and then set
\[f=\dots f_{1-x}^{(a_4)}f_x^{(a_3)}f_{1-x}^{(a_2)}f_x^{(a_1)}\in\calu.\]
(For $j$ sufficiently large we have $\mu^j=\la$, and so $a_j=0$ for large $j$, so this definition makes sense.)  Our objective is to compute the coefficient of $s(\la)$ in $fs(\mu)$.  To do this, we construct a $\la$-tableau $T$ by filling each node of $[\mu]$ with a $0$, and then for $j\gs1$, filling each node of $[\mu^j]\setminus[\mu^{j-1}]$ with a $j$.  Given a node $(k,h)\in[\la]$, let $j=T(k,h)$, and set
\begin{align*}
N(k,h) = &\left|\left\{m<k\ \left|\ T(m,\la_m)<j,\ T(m,\la_m)\nequiv j\pmod2\right.\right\}\right|\\
-&\left|\left\{m<k\ \left|\ T(m,\la_m)<j,\ T(m,\la_m)\equiv j\pmod2\right.\right\}\right|.
\end{align*}
Finally, set $N=\sum_{(k,h)\in[\la]}N(k,h)$.  Now we we have the following.

\begin{lemma}\label{flem}
With the above definitions, we have
\[\lan fs(\mu),s(\la)\ran = v^N.\]
\end{lemma}

\begin{pf}
The hypothesis on $\mu$ means that for any $i$, the nodes $(i,\mu_i+1),(i,\mu_i+2),\dots,(i,\la_i)$ are filled with the integers $1,2,\dots,\la_i-\mu_i$ in $T$.  In particular, for each $j\gs1$, the nodes $(k,h)$ with $T(k,h)=j$ all have residue $x+j\pmod2$.  Now the lemma is straightforward to prove, given the above formula for the actions of $f^{(a)}_0,f^{(a)}_1$ on $\calf$.
\end{pf}

\begin{eg}
Set $\la=(7^2,5^2,4)$, $\mu=(7,6,3,2,1)$, $x=1$.  Then we have
\begin{align*}
\mu^0 &=\mu,\\
\mu^1 &=(7^2,4,3,2),\\
\mu^2 &=(7^2,5,4,3),\\
\mu^j &=\la\text{ for $j\gs3$},
\end{align*}
and
\begin{align*}
T =\ &\raisebox{-27pt}{\young(0000000,0000001,00012,00123,0123)}.\\
\intertext{The values of $N(k,l)$ are given by}
&\raisebox{-27pt}{\young(0000000,0000001,00010,00101,0101)},
\end{align*}
so that
\[\left\lan f_1^{(2)}f_0^{(3)}f_1^{(4)}s\big((7,6,3,2,1)\big),s\left((7^2,5^2,4)\right)\right\ran = v^6.\]
\end{eg}

We now use the last two results to prove the following proposition, which covers the bulk of cases in our main theorem.

\begin{propn}\label{fockprop}
Suppose $\la$ is a partition satisfying the hypothesis of Theorem \ref{main}.  Let $l=\la'_1$, and suppose that for $i=1,\dots,l$ we have $\la_i\gs l-i+2$. Suppose also that for $1\ls k\ls \la_1$ the ladder $\call_k(\la)$ is connected. Then $\spe\la $ is reducible.
\end{propn}

\begin{pf}
By the results of \S\ref{adjsec}, we may assume that $\bbf=\bbq$.  The hypothesis of Theorem \ref{main} is that some ladder $\call_m(\la)$ is disconnected; take the smallest such $m$, and choose $i$ such that $(i-1,m+2-i)\in[\la]\notni(i,m+1-i)$ and $(j,m+1-j)\in[\la]$ for some $j>i$.  Now define
\begin{align*}
\mu &= (l+1,l,l-1,\dots,3,2),\\
\hat\mu &= (l+1,l,l-1,\dots,l-i+3,l-i,l-i-1,\dots,2,1).
\end{align*}
Setting $x=l+1\pmod2$, we find that $\mu$ and $\hat\mu$ both satisfy ($\ast$); we define the operator $f$ and the tableau $T$ corresponding to $(\la,\mu,x)$ as above, we define $\hat f$ and $\hat T$ corresponding to $(\la,\hat\mu,x)$ in the same way.  Our aim is to show that the hypotheses of Lemma \ref{mctsl} are satisfied, with
\[X = fG(\mu),\qquad Y=\hat fG(\hat\mu).\]
Certainly $X$ and $Y$ are bar-invariant elements of the Fock space.  We now claim that we can ignore all terms in $G(\mu)$, $G(\hat\mu)$ except the leading terms, i.e.
\[\lan X,s(\la)\ran = \lan fs(\mu),s(\la)\ran,\qquad \lan Y,s(\la)\ran = \lan \hat fs(\hat\mu),s(\la)\ran.\]
Note that $\mu$ is an alternating partition with core $(l-1,l-2,\dots,1)$; so by Corollary \ref{altcory} and Theorem \ref{arikithm}, any $\nu\neq\mu$ which gives a non-zero term $d_{\nu\mu}(v)s(\nu)$ in $G(\mu)$ satisfies $\nu_{l+1}>0$.  But this means that $[\nu]\nsubseteq[\la]$, which obviously implies that $\lan fs(\nu),s(\la)\ran=0$.  The same argument applies to $G(\hat\mu)$, and so we can concentrate on $fs(\mu)$ and $fs(\hat\mu)$.  If we define the integer $N$ corresponding to $T$ as above, and define $\hat N$ from $\hat T$ analogously, then by Lemma \ref{flem} we have
\[\lan fs(\mu),s(\la)\ran=v^N,\qquad\lan \hat fs(\hat\mu),s(\la)\ran=v^{\hat N},\]
and it remains to prove the purely combinatorial statement that $\hat N\neq N$.

In fact, we shall estimate $\hat N-N$, and show that it is strictly positive.  To do this, we compare $N(k,h)$ with $\hat N(k,h)$ for the various nodes $(k,h)\in[\la]$.  It will help to introduce some notation: for any $j\gs0$, we let $a_j$ be the number of $k\in\{1,\dots,i-1\}$ such that $T(k,\la_k)=j$; that is, the number of rows of $T$ above row $i$ ending in $\rh{\young(j)}$.  We also define $b_j$ to be the number of nodes $(k,h)$ with $k\gs i$ such that $T(k,h)=j$.

First we note that $T$ and $\hat T$ agree on rows $1,\dots,i-1$, so certainly for any node $(k,h)$ with $k<i$ we have $\hat N(k,h)=N(k,h)$.  So we concentrate on nodes in rows $i$ and below.  Consider first the nodes lying in $[\mu]\setminus[\hat\mu]$.  There are two of these in each row from $i$ to $l$; these are labelled $\rh{\young(00)}$ in $T$, and $\rh{\young(12)}$ in $\hat T$.  By the definitions of $N$ and $\hat N$, each such pair of nodes contributes nothing to $N$, and contributes $a_1$ to $\hat N$ (the node $\rh{\young(1)}$ contributes $a_0$, while the node $\rh{\young(2)}$ contributes $-a_0+a_1$).

Next let $(k,h)$ be a node of $[\la]\setminus[\mu]$ with $k\geq i$, and let $j=T(k,h)=k+h-l-2$.  Then $\hat T(k,h)=j+2$, and using the definitions we can compute
\[\hat N(k,h)-N(k,h) = a_{j+1}-a_j\]
(note that the rows from $i$ to $k-1$ make no contribution).  Summing, we find that
\begin{align*}
\hat N-N &=(l-i+1)a_1 + (a_2-a_1)b_1 + (a_3-a_2)b_2+\dots\\
&= (l-i+1-b_1)a_1+(b_1-b_2)a_2+(b_2-b_3)a_3+\dots.
\end{align*}
We claim that each term of the latter sum is non-negative, and that some term is positive.  Certainly $b_1\ls l-i+1$, since there can be at most one node labelled $1$ in any of rows $i,\dots,l$.  Also, $b_j\ls b_{j-1}$ for $j\gs2$, since each node labelled $j$ must have a node labelled $j-1$ immediately to its left.  So each term of the sum is non-negative.  Now let $j=T(i-1,m-i+2)$; that is, $j=m-l-1$.  If $j\geq 2$, we claim that $(b_{j-1}-b_j)a_j>0$; otherwise we claim $(l-i+1-b_1)a_1>0$.  Suppose $j\geq 2$.  To see that the factor $b_{j-1}-b_j$ is positive, note that there is a node labelled $j-1$ at the end of row $i$, but no node labelled $j$ in this row; in any subsequent row, if there is a node labelled $j$ then there is a node labelled $j-1$ immediately to its left.  Similarly if $j=1$, there is no node in row $i$ labelled 1, and so $l-i+1>b_1$.  

Now we show that $a_j>0$.  The last hypothesis of the proposition implies that ladder $m$ does not meet row $1$, so that $T(1,\la_1)<j$; on the other hand, there is a node labelled $j$ in row $i-1$, so $T(i-1,\la_{i-1})\gs j$.  For any $1<k\ls i-1$ it is easy to see that $T(k,\la_k)\ls T(k-1,\la_{k-1})+1$ so every value from $T(1,\la_1)$ to $T(i-1,\la_{i-1})$ occurs as some $T(k,\la_k)$ for $1\ls k\ls i-1$.  In particular, the value $j$ occurs, and we are done.
\end{pf}

\sect{Proof of Theorem \ref{main}}\label{pfsec}

\begin{propn}\label{twoends}
Suppose $\la$ is a partition and $l\gs1$.  Suppose that $(1,l)$ and $(l,1)$ both lie in $\call_l(\la)$, but that $\call_l(\la)$ is disconnected.  Then $\spe\la$ is reducible.
\end{propn}

\begin{pf}
Suppose $\la$ has weight $w$ and core $(r,r-1,\dots,1)$; we proceed by induction on $w$, and for fixed $w$ by reverse induction on $r$.  The starting case for this induction is where $r\gs w-1$, so that $\la$ lies in a Rouquier block.  Since $\la$ is certainly neither $2$-regular nor $2$-restricted, Corollary \ref{rouqcory} gives the result in this case.

For the case where $r<w-1$, recall the definition of the partition $\la^{+i}$ for $i\in\{0,1\}$ from \S\ref{213sec}.  By Lemma \ref{wtdown}, we have $w(\la^{+i})\ls w(\la)$, and obviously if $w(\la^{+i})=w(\la)$ but $\la^{+i}\neq\la$ then the core of $\la^{+i}$ is larger than the core of $\la$.  So to complete the inductive step is suffices to show that for either $i=0$ or $i=1$ we have $\la^{+i}\neq\la$ with $\la^{+i}$ also satisfying the hypotheses of the lemma.

Let $j=l\pmod2$.  Consider two cases.
\begin{itemize}
\item
Suppose $\la$ has at least one addable node of residue $j$.  Then we have $\la^{+j}\neq\la$, and (since all of the nodes in $\call_l$ have residue $1-j$) we have $\call_l(\la^{+j})=\call_l(\la)$, so $\la^{+j}$ satisfies the hypotheses of the lemma.
\item
Alternatively, suppose $\la$ has no addable nodes of residue $j$.  Then the nodes $(1,l+1)$ and $(l+1,1)$ must lie in $[\la]$ (since these nodes have residue $j$, and if either of them were not contained in $[\la]$ then it would be an addable node).  On the other hand, $[\la]$ cannot contain all the nodes in $\call_{l+1}$ since it does not contain all the nodes in $\call_l$; so ladder $\call_{l+1}$ is disconnected.  Now we can replace $l$ with $l+1$ (and $j$ with $1-j$), and appeal to the previous case.
\end{itemize}
\end{pf}

\begin{propn}\label{oneend}
Suppose $\la$ is a partition and $l\gs1$.  Suppose the ladder $\call_l(\la)$ is disconnected, and that $(1,l)\in[\la]\notni(l,1)$.  Then $\spe\la$ is reducible.
\end{propn}

\begin{pf}
Divide $\call_l(\la)$ into `segments' of consecutive nodes; the condition that $\call_l(\la)$ is disconnected is precisely the statement that there are at least two segments.  Let $s$ be the length of the shortest segment other than the segment containing the node $(1,l)$, and proceed by induction on $s$.  Let $j=l+1\pmod2$ be the common residue of the nodes in $\call_l$, and let $\la^+=\la^{+(1-j)}$ be the partition defined in \S\ref{213sec}.

Suppose $s=1$.  This means that there is some $i\in\{3,\dots,l-1\}$ such that $(i,l+1-i)$ is a node of $\la$ but neither $(i-1,l+2-i)$ nor $(i+1,l-i)$ is.  In particular, this implies that $(i,l+1-i)$ is a removable node of $\la$.  We claim that this node is also a removable node of $\la^+$: since neither of the nodes $(i-1,l+2-i)$, $(i+1,l-i)$ lies in $[\la]$ by assumption, neither of the nodes $(i,l+2-i)$, $(i+1,l+1-i)$ can lie in $[\la^{+(1-j)}]$; so $(i,l+1-i)$ is a removable node of $\la^+$, as claimed.  Now consider the addable node $\fka=(1,\la^+_1+1)$ of $\la^+$.  Since $\la^+$ cannot have addable nodes of residue $1-j$, $\fka$ must have residue $j$.  Moreover, since $(1,l)\in[\la]$, $\fka$ lies in ladder $\call_m$ for some $m>l$.  Hence by Proposition \ref{ladhom} $\spe{\la^+}$ is reducible, and so by Lemma \ref{213dual} $\spe\la$ is reducible.

Now we consider the inductive step; suppose $s>1$, and that
\[\{(i,l+1-i),(i+1,l-i),\dots,(i+s-1,l+2-i-s)\}\]
is a segment of length $s$, with $3\ls i\ls l-s$.  That is, the nodes listed above are nodes of $\la$, but the nodes $(i-1,l+2-i)$, $(i+s,l+1-i-s)$ are not.  This means that the nodes
\[(i+1,l+1-i),(i+2,l-i),\dots,(i+s-1,l+3-i-s)\]
are either nodes or addable nodes of $\la$, and hence (since they have residue $1-j$) are nodes of $\la^+$.  On the other hand, neither of the nodes $(i,l+2-i),(i+s,l+2-i-s)$ is a node or an addable node of $\la$, so neither of these nodes is a node of $\la^+$.  So $\call_{l+1}(\la^+)$ includes a segment of length $s-1$; it also contains the node $(1,l+1)$ but not the node $(l+1,1)$, and so we may apply the inductive hypothesis, replacing $\la$ with $\la^+$ and $l$ with $l+1$, to deduce that $\spe{\la^+}$ is reducible.  Now we can apply Lemma \ref{213dual}.
\end{pf}

\begin{eg}
Let $\la=(5,3^2,2)$.  Then $\la$ satisfies the hypotheses of Proposition \ref{oneend}, with $l=5$.  We have $s=2$ and $j=0$, and we examine the partition $\la^+=\la^{+1}=(6,3^3)$.  This partition also satisfies the hypotheses of Proposition \ref{oneend}, with $l=6$ and $s=1$.  We construct the partition $\la^{++}=(\la^+)^{+0} = (7,4,3^2,1)$; this satisfies the hypotheses of Proposition \ref{ladhom}, since it has an addable node $(1,8)\in\call_8$ and a removable node $(4,3)\in\call_6$.  So $\spe{\la^{++}}$ is reducible, and hence so is $\spe{\la^+}$, and hence so is $\spe\la$.  The Young diagrams of the partitions, with the residues of their nodes marked, are given below.
\[\begin{array}{c@{\qquad\qquad}c@{\qquad\qquad}c}
\la&\la^+&\la^{++}\\[6pt]
\young(01010,101,010,10)&\young(010101,101,010,101)&\raisebox{-12pt}{\young(0101010,1010,010,101,0)}
\end{array}\]
\end{eg}

\begin{propn}\label{square}
Suppose $\la$ satisfies the hypothesis of Theorem \ref{main}, that $\la_1=\la'_1$, and that $\la_i\gs \la_1+1-i$ for $i=1,\dots,\la_1$.  Then $\spe\la$ is reducible.
\end{propn}

\begin{pf}
Using Proposition \ref{mattshom} and Corollary \ref{homcory}, our task is to show that we can find partitions $\nu,\xi$ such that $\xi_{\nu'_1-1}\gs \nu'_1$, $\la_i = \xi_i+2\nu_i$ for all $i$, and $\la\reg\ndom\mu$, where $\mu$ is the partition given by $\mu_i = \xi'_i+2\nu_i$ for all $i$.

Let $m$ be minimal such that $\call_m(\la)$ is disconnected, and let $i,k$ be such that $k\gs i+2$,  $(i,m+1-i),(k,m+1-k)$ lie in $[\la]$ and none of the nodes $(i+1,m-i),\dots,(k-1,m+2-k)$ lies in $[\la]$.  We would like to assume that $i+1\ls m-i$, i.e.\ the node $(i+1,m-i)$ lies on or above the main diagonal of the Young diagram.  If this is not the case, then we can replace $\la$ with $\la'$ (appealing to Corollary \ref{conjirr}), and replace $(m,i,k)$ with $(m,\bar\imath,\bar k)$, where $\bar\imath=m+1-k$, $\bar k=m+1-i$; it is then easy to check that $\bar\imath+1\ls m-\bar\imath$.

So we shall assume that $i+1\ls m-i$.  Since $(i,m+1-i)\in[\la]\notni(i+1,m-i)$, we have $\la_i-\la_{i+1}\gs2$.  So we if we define $\nu=(1^i)$ and $\xi=(\la_1-2,\la_2-2,\dots,\la_i-2,\la_{i+1},\la_{i+2},\dots)$, then $\xi$ is a partition.  Furthermore, we have 
\[\xi_{\nu'_1-1} = \xi_{i-1} = \la_{i-1}-2\gs \la_i-2\gs m-1-i\gs i=\nu'_1,\]
and it remains to show that $\la\reg\ndom\mu$.  In fact, we shall show that $\mu'_1<(\la\reg)'_1$, which is certainly good enough.  By assumption $\call_{\la_1}(\la)=\call_{\la_1}$, and hence $\call_{\la_1}(\la\reg)=\call_{\la_1}$, and this means that $(\la\reg)'_1=\la_1$.  On the other hand,
\[\mu'_1 = \max\{\xi_1,\nu'_1\} = \max\{\la_1-2,i\}<\la_1,\]
and we are done.
\end{pf}

\begin{pfof}{Theorem \ref{main}}
Suppose $\la$ satisfies the hypothesis of Theorem \ref{main}.  By replacing $\la$ with $\la'$ if necessary and appealing to Corollary \ref{conjirr}, we may assume that $\la_1\gs\la'_1$.

Let $m$ be minimal such that ladder $\call_m(\la)$ is disconnected.  If $m\ls\la_1$, then $\call_m(\la)$ meets the top row of $[\la]$ (i.e.\ $(1,m)\in[\la]$), and so we may appeal to Proposition \ref{twoends} or Proposition \ref{oneend}.  So we can assume that $m>\la_1$, and in particular no disconnected ladder of $\la$ meets the top row.

If $\la_1=\la'_1$, then we may appeal to Proposition \ref{square}, so instead we suppose that $\la_1>\la'_1$.  Let $l=\la'_1$, and suppose that $\la_l\gs2$.  Then we have $(l,2)\in[\la]$.  Since $\la_1>l$, we also have $(1,l+1)\in[\la]$, and now the assumption that $\call_{l+1}(\la)$ is connected means that the nodes $(2,l),(3,l-1),\dots,(l-1,3)$ all lie in $\la$.  So $\la_i\gs l-i+2$ for $i=1,\dots,l$, and we may appeal to Proposition \ref{fockprop}.

We are left with the case where $\la_l=1$.  In this case, $\la$ has a removable node $(l,1)$, of residue $j=l+1\pmod2$.  We now consider two cases.
\begin{itemize}
\item
Suppose the nodes in ladder $\call_m$ have residue $j$.  We claim that there is an addable node of $\la$ in this ladder.  Indeed, let $i,k$ be such that $k\gs i+2$, $(i,m+1-i),(k,m+1-k)$ lie in $[\la]$ and none of the nodes $(i+1,m-i),\dots,(k-1,m+2-k)$ lies in $[\la]$.  Then $(i,m-i)$ and $(k-1,m+1-k)$ lie in $[\la]$, and since $\call_m(\la)$ is the first disconnected ladder of $\la$, the nodes $(i+1,m-1-i),\dots,(k-2,m+2-k)$ lie in $[\la]$.  So $\la$ has an addable node $(i+1,m-i)\in\call_m$.  Now we may appeal to Proposition \ref{ladhom}.
\item
Alternatively, suppose the nodes in ladder $\call_m$ have residue $1-j$, and consider the partition $\la^{-j}$.  This is strictly smaller than $\la$ since it does not contain the node $(l,1)$, but also has a disconnected ladder, i.e.\ $\call_m(\la^{-j})=\call_m(\la)$.  So by induction on $|\la|$ and Lemma \ref{213}, $\spe\la$ is reducible.
\end{itemize}
\end{pfof}

\end{document}